\documentclass[11pt]{gtart}
\usepackage{amsmath, amssymb, epsfig, graphs, psfrag}

\newtheorem{thm}{Theorem}[section]

\newtheorem{cor}[thm]{Corollary}
\newtheorem{lem}[thm]{Lemma}
\newtheorem{prop}[thm]{Proposition}

\theoremstyle{definition}

\newtheorem{rem}[thm]{Remark}

\numberwithin{equation}{section}

\newcommand{\Ga}{\Gamma}

\newcommand{\de}{\delta}

\newcommand{\om}{\omega}

\newcommand{\si}{\sigma}
\newcommand{\Si}{\Sigma}

\newcommand{\x}{\times}
\newcommand{\s}{\mathbf s}

\renewcommand{\t}{\mathbf t}

\newcommand{\Z}{\mathbb Z}
\newcommand{\N}{\mathbb N}
\newcommand{\Q}{\mathbb Q}
\newcommand{\R}{\mathbb R}

\newcommand{\CP}{{\mathbb C}{\mathbb P}}
\newcommand{\cpkk}{{\overline {{\mathbb C}{\mathbb P}^2}}}

\newcommand{\del}{\partial}

\newcommand{\lra}{\longrightarrow}
\newcommand{\hf}{{{\widehat {HF}}}}

\begin{document}

\title{Ozsv\'ath--Szab\'o invariants and\\ 
tight contact three--manifolds, III}

\author{Paolo Lisca}
\address{Dipartimento di Matematica\\
Universit\`a di Pisa \\I-56127 Pisa, ITALY} 
\email{lisca@dm.unipi.it}

\author{Andr\'{a}s I. Stipsicz}
\address{R\'enyi Institute of Mathematics\\
Hungarian Academy of Sciences\\
H-1053 Budapest\\ 
Re\'altanoda utca 13--15, Hungary and \\
Institute for Advanced Study, Princeton, NJ}
\email{stipsicz@math-inst.hu}

\begin{abstract}
We characterize $L$--spaces which are Seifert fibered over the
2--sphere in terms of taut foliations, transverse foliations and
transverse contact structures. We give a sufficient condition for
certain contact Seifert fibered 3--manifolds with $e_0=-1$ to have
nonzero contact Ozsv\'ath--Szab\'o invariants. This yields an
algorithm for deciding whether a given small Seifert fibered
$L$--space carries a contact structure with nonvanishing contact
Ozsv\'ath--Szab\'o invariant. As an application, we prove the
existence of tight contact structures on some 3--manifolds obtained by
integral surgery along a positive torus knot in the
3--sphere. Finally, we prove planarity of every contact structure on
small Seifert fibered $L$--spaces with $e_0\geq -1$, and we discuss
some consequences.
\end{abstract}

\maketitle

\section{Introduction}
The Ozsv\'ath--Szab\'o homology groups of a closed, oriented
3--manifold $Y$~\cite{OSzF1, OSzF2} capture important topological
information about $Y$. For example, by~\cite[Theorem~1.1]{OSzgen} the
Thurston semi--norm is determined by the evaluation of the first Chern
classes of spin$^c$ structures with nontrivial Ozsv\'ath--Szab\'o
homology groups. For rational homology spheres, however, the Thurston
norm is trivial, while the Ozsv\'ath--Szab\'o homology groups are
special: the group $\hf (Y ,\t )$ has odd rank for each spin$^c$
structure $\t\in Spin ^c(Y)$. A rational homology sphere $Y$ shows the
simplest possible Heegaard Floer--theoretic behaviour if for every 
spin$^c$ structure $\t\in Spin ^c(Y)$ the group $\hf (Y, \t )$ with
integral coefficients is isomorphic to $\Z$, in which case $Y$ is
called an \emph{$L$--space}. In the present paper we shall always use
$\Z/2\Z$--coefficients. In this case, if $Y$ is an $L$--space then
$\hf(Y,\t)\cong\Z/2\Z$ for every $\t\in Spin^c(Y)$.

If $Y$ is a 3--manifold with elliptic geometry then $Y$ is an
$L$--space~\cite[Proposition~2.3]{OSzlens}. In particular, lens spaces
are $L$--spaces. Seifert fibered $L$--spaces can be characterized
combinatorially using the results of~\cite{Ne, OSzplum}.

More generally, by~\cite{Ne} a negative definite plumbed 3--manifold
is an $L$--space if it is the link of a rational surface
singularity. In~\cite[Section~5, Question 11]{OSzsurvey} it was asked
whether there is a topological characterization of $L$--spaces. In the
first part of this paper we study the topological significance of
being an $L$--space for Seifert fibered rational homology 3--spheres.

It is proved in~\cite[Theorem~1.4]{OSzgen} that an $L$--space admits
no taut foliation. Our first result shows that the converse to this
statement holds for Seifert fibered rational homology
3--spheres.  In addition, we find that it is equivalent to the
nonexistence of a transverse foliation (i.e., a foliation transverse
to the fibers of the Seifert fibration), as well as to the
nonexistence of positive, transverse contact structures for at least
one orientation.

\begin{figure}[ht]
\begin{center}
\psfrag{r1}{$-\frac1{r_1}$}
\psfrag{r2}{$-\frac1{r_2}$}
\psfrag{rk}{$-\frac1{r_k}$}
\psfrag{e0}{$e_0$}
\includegraphics[height=2.5cm]{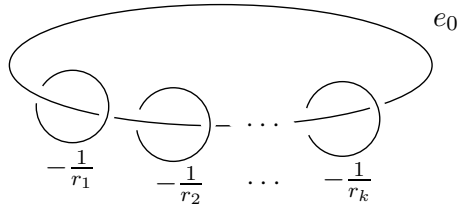}
\end{center}
\caption{\quad Surgery diagram for the Seifert fibered 3--manifold
$M(e_0; r_1, r_2,\ldots, r_k)$}
\label{f:seifert}
\end{figure}

Let $M$ be a closed, oriented rational homology 3--sphere which
carries a Seifert fibration over $S^2$. Then, $M$ is orientation
preserving diffeomorphic to the oriented 3--manifold $M(e_0; r_1,
r_2,\ldots, r_k)$ with surgery presentation given by
Figure~\ref{f:seifert}, where $e_0\in \Z$, $r_i \in (0,1)\cap \Q$ and
$r_1\geq r_2\geq\cdots\geq r_k$. Moreover, the surgery picture shows
that $M(e_0;r_1,\ldots, r_k)$ carries a natural Seifert fibration over
$S^2$, which coincides (up to isotopy) with the image of the fibration
on $M$ under the above diffeomorphism. Througout the paper, we 
always implicitely refer to this natural fibration on
$M(e_0;r_1,\ldots,r_k)$.

Define 
\[
e(M):=e_0(M)+r_1+r_2+\cdots+r_k.  
\]
Notice that
\[
-M(e_0,r_1,r_2,\ldots,r_k)=M(-e_0-k; 1-r_1,1-r_2,\ldots,1-r_k),
\]
therefore we have $e(-M)=-e(M)$. It is known that $M$ is a rational
homology sphere if and only if $e(M)\neq 0$.

\begin{thm}\label{t:equiv}
Let us suppose that $M$ is an oriented rational homology 3--sphere
which is Seifert fibered over $S^2$. Then, the following statements are
equivalent:
\begin{enumerate}
\item $M$ is an $L$--space
\item Either $M$ or $-M$ carries no positive, transverse contact
structures
\item $M$ carries no transverse foliations
\item $M$ carries no taut foliations.
\end{enumerate}
\end{thm}

Recall that 
\[
\hf (M,\t)\cong \hf (-M,\t) 
\]
for each $\t\in Spin^c(M)$. It was proved in~\cite{LM} that an
oriented Seifert fibered rational homology 3--sphere $M=M(e_0,r_1,r_2,
\ldots ,r_k)$ with $r_1\geq r_2\geq \cdots \geq r_k$ admits no
positive transverse contact structure if and only if
\begin{itemize}
\item $e_0(M)\geq 0$, or
\item $e_0(M)=-1$ and there are no relatively prime integers
$m>a$ such that
\[
mr_1<a<m(1-r_2)\quad\text{and}\quad  mr_i<1 \quad i=3, \ldots , k.
\]
\end{itemize}

Combined with this result, Theorem~\ref{t:equiv} gives a simple
characterization of $L$--spaces among rational homology
3--spheres of the form $M(e_0;r_1,\ldots,r_k)$.

In the second part of the paper we give a sufficient condition for
certain contact Seifert fibered 3--manifolds with $e_0=-1$ to have
nonzero contact Ozsv\'ath--Szab\'o invariants. Combining this result
with Theorem~\ref{t:equiv} and with a result of Paolo Ghiggini
(Theorem~\ref{t:ghig} below) yields an algorithm for deciding whether
an $L$--space of the form $M(e_0;r_1,r_2,r_3)$ carries a contact
structure with nonvanishing contact Ozsv\'ath--Szab\'o invariant.

\begin{figure}[ht]
\begin{center}
\psfrag{+1}{\small $+1$}
\psfrag{r1}{\small $-\frac1{r_1}$}
\psfrag{r2}{\small $-\frac1{r_2}$}
\psfrag{rk}{\small $-\frac1{r_k}$}
\epsfig{file=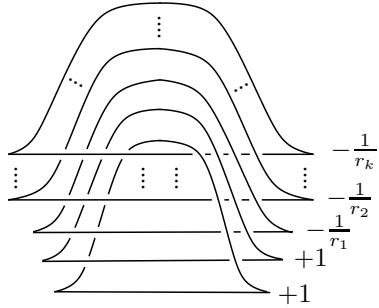, height=4cm}
\end{center}
\caption{Contact structures on $M(-1;r_1,r_2,\ldots,r_k)$}
\label{f:special}
\end{figure}

Consider the contact surgery diagram of Figure~\ref{f:special}.  The
diagram should be interpreted as representing all possible contact
$(-\frac{1}{r_i})$--surgeries in case these are not unique (see
Section~\ref{s:sec} for details).
With this convention, Figure~\ref{f:special} represents a family of
contact structures on the manifold $M=M(-1;r_1, \ldots ,r_k)$. Given an
oriented contact structure $\xi$ on $M$, let $\t _{\xi }$ denote the
spin$^c$ structure induced by $\xi$ and $d_3(\xi )\in\Q$ the 3--dimensional
invariant determined by the homotopy class of the oriented 2--plane
field $\xi$~\cite{Go}. Finally, denote by $d(M,\t)\in\Q$ the 
$d$--invariant of $(M, \t)$ (cf. \cite{abs}).
\begin{thm}\label{t:suff}
Let $(M,\xi )$ be any of the contact, Seifert fibered
3--manifolds given by Figure~\ref{f:special}. Suppose that 
\[
e(M)=r_1+r_2+\cdots+r_k-1>0, 
\]
and 
\[
d_3(\xi)=d(M, \t _{\xi}).
\]
Then the contact Ozsv\'ath--Szab\'o invariant $c(M,\xi )\in \hf (-M,
\t _{\xi })$ is nonzero.
\end{thm}

Theorem~\ref{t:suff} applies to some Seifert fibered 3--manifolds
which are not $L$--spaces (such as
e.g.~$M(-1;\frac5{12},\frac13,\frac13)$), but it is particularly
useful when $k=3$ and the 3--manifold $M$ is an $L$--space. To explain why, we
need a small digression.

Suppose that $(M,\xi)$ is a contact, Seifert fibered
$3$--manifold. Then, a Legendrian knot in $M$ smoothly isotopic to a
regular fiber admits two framings: one coming from the fibration and
another one coming from the contact structure $\xi$. The difference
between the contact framing and the fibration framing is the
\emph{twisting number} of the Legendrian curve. We say that $\xi$ has
\emph{maximal twisting equal to zero} if there is a Legendrian knot
$L$ isotopic to a regular fiber such that $L$ has twisting number
zero. According to a result of Ghiggini, if a rational homology
3--sphere $M=M(e_0;r_1,r_2,r_3)$ admits no transverse contact
structure then each tight contact structure on $M$ has maximal
twisting equal to zero:

\begin{thm}[Ghiggini, \cite{Gh}]\label{t:ghig}
If $M=M(e_0;r_1,r_2,r_3)$ is a rational homology 3--sphere and admits
a tight contact structure with negative twisting, then $M$ admits a
transverse contact structure.\qed
\end{thm}

Thus, combining Theorems~\ref{t:equiv} and~\ref{t:ghig} we see that
\emph{if $M=M(e_0;r_1,r_2,r_3)$ is an $L$--space, then each tight
contact structure $\xi$ on $M$ has maximal twisting equal to zero}. On
the other hand, in Section~\ref{s:six} (see Proposition~\ref{p:model})
we show that every tight contact structure with maximal twisting equal
to zero on a 3--manifold of the form $M(e_0;r_1,r_2,r_3)$ is given by
one of the diagrams represented by Figure~\ref{f:special}. Moreover,
the quantities $d_3(\xi)$ and $d(M,\t_{\xi})$ can be computed
algorithmically from a contact surgery presentation. Therefore
Theorem~\ref{t:suff} yields an algorithm for deciding whether an
$L$--space of the form $M(e_0;r_1,r_2,r_3)$ carries a contact
structure with nonvanishing contact Ozsv\'ath--Szab\'o invariants.
(Recall that contact structures with nonvanishing contact
Ozsv\'ath--Szab\'o invariant are tight.)  In fact, one can view
Theorem~\ref{t:suff} as a useful tool to attack the classification
problem for tight contact structures on small Seifert fibered
$L$--spaces with $e_0=-1$. (Tight contact structures on small Seifert
fibered 3--manifolds with $e_0\neq -1, -2$ have been classified
\cite{GLS, Wu}, but the classification for $e_0=-1,-2$ is expected to
be considerably harder because of the presence of nonfillable
structures~\cite{EH1, LS1} as well as manifolds with no tight contact
structures~\cite{EH2, LS2}.)

As an application of Theorem~\ref{t:suff} we prove the existence
of tight contact structures on some of the 3--manifolds obtained by
integral surgery along a positive torus knot in $S^3$. To put this
result in perspective, recall that we showed in~\cite{LSgt} that for
the positive $(p,q)$--torus knot $T_{p,q}\subset S^3$ the 3--manifold
$S^3_r (T_{p,q})$ obtained by topological $r$--surgery along $T_{p,q}$
carries tight contact structures for every rational $r$ except
possibly when $r=pq-p-q$. In this latter case the idea used
in~\cite{LSgt} does not work, since it would require the application
of contact 0--surgery, which automatically leads to overtwisted
structures. In~\cite{LS2} we analyzed the existence of tight contact
structures on the 3--manifold $S^3_r(T_{p,q})$ with $r$ equal to the
`critical' surgery coefficient $pq-p-q$. We showed that for
$(p,q)=(2,2n+1)$ the corresponding manifolds do not admit any tight
structures (extending an earlier result of Etnyre and Honda
\cite{EH2}, who proved the same result for $(p,q)=(2,3)$). On the
other hand, in the same paper~\cite{LS2} we verified that
$S^3_{p^2n-pn-1}(T_{p,pn+1})$ does carry tight contact structures for
every $n\geq 1$ and $p\geq 3$ \emph{odd}. The parity assumption on $p$
played a crucial role in the argument, since in this case the
candidate contact structure induced a spin structure.  This allowed us
to use a certain $\Z /2\Z $ symmetry (the so--called $J$--map) built
in Heegaard Floer theory to prove the nonvanishing of the contact
invariant. A little computer search using the algorithm outlined above
showed that for $p$ even one cannot always find a tight contact
structure with nonzero contact Ozsv\'ath--Szab\'o invariant inducing a
spin structure. For example, there is no such structure on
$S^3_{11}(T_{4,5})$.  This implies that when $p$ is even, a different
approach should be used. In Section~\ref{s:critical} we use
Theorem~\ref{t:suff} to prove:

\begin{thm}\label{t:existence}
For each $p>2$, the manifold $S^3_{p^2-p-1}(T_{p,p+1})$ carries tight
contact structures.
\end{thm}

Using the fact that every tight contact structure with zero maximal twisting 
on a 3--manifold of the form $M(e_0;r_1,r_2,r_3)$ is given by a surgery 
diagram as in Figure~\ref{f:special} (for $k=3$), we can prove: 

\begin{thm}\label{t:planob}
A tight contact structure $\xi$ with zero maximal twisting on a small
Seifert fibered 3--manifold $M(e_0; r_1,r_2,r_3)$ is compatible with a
planar open book decomposition.
\end{thm}

An immediate corollary of Theorem~\ref{t:planob} is the following 
characterization of small Seifert fibered $L$--spaces:

\begin{cor}\label{c:planar=L}
A rational homology 3--sphere $M=M(e_0;r_1,r_2,r_3)$ is an $L$--space
if and only if on either $M$ or $-M$ every contact structure is
planar.
\end{cor}

Another corollary of Theorem~\ref{t:planob} is the following:

\begin{cor}\label{c:allplanar}
If $M=M(e_0;r_1,r_2,r_3)$ with $e_0\geq -1$ and $M$ is an $L$--space,
then each contact structure on $M$ is planar.
\end{cor}

According to a recent result of Abbas, Cieliebak and Hofer~\cite{ACH}
this fact implies

\begin{cor}\label{c:weinstein}
If $M=M(e_0;r_1,r_2,r_3)$ with $e_0\geq -1$ and $M$ is an $L$--space,
then each contact structure $\xi$ on $M$ satisfies the Weistein
conjecture, that is, every Reeb vector field for $(M,\xi)$ admits a
periodic orbit.
\qed\end{cor}

The paper is organized as follows. In Section~\ref{s:sec} we collect
some background material. In Sections~\ref{s:firstproof},
\ref{s:tight} and~\ref{s:critical} we prove respectively
Theorems~\ref{t:equiv}, \ref{t:suff} and~\ref{t:existence}.
Section~\ref{s:six} is devoted to examining compatible open books and
establishing Theorem~\ref{t:planob} and its corollaries.

{\bf Acknowledgements.}  The authors are grateful to Stefan Friedl for
suggesting the neat argument used in the proof of
Proposition~\ref{p:alex}. The first author was partially supported by
MURST. The second author was partially supported by OTKA T049449.

\section{Preliminaries}\label{s:sec}

\sh{Transverse contact structures and transverse foliations} 

Here we collect some known results on transverse contact structures
and foliations on the  Seifert fibered 3--manifold 
\[
M=M(e_0; r_1, r_2,\ldots, r_k),\quad
r_1\geq r_2\geq\cdots\geq r_k 
\]
defined by Figure~\ref{f:seifert}. 
Consider
\[
\Ga(M):=(r_1,r_2,\ldots,r_k)\in(\Q\cap(0,1))^k.
\]
We say that $\Ga(M)$ is~\emph{realizable} if there exist coprime
integers $m>a>0$ such that
\[
r_1 < \frac am, \quad
r_2 < \frac{m-a}m, \quad\text{and}\quad
r_3,\ldots,r_k < \frac 1m.
\]
In Section~\ref{s:firstproof} we shall use the following theorem,
obtained by combining results of Jenkins--Neumann and Naimi, and stated
here in our present notation:

\begin{thm}[\cite{JN, Na}]\label{t:quoteJN}
Let $M=M(e_0;r_1,r_2,\ldots,r_k)$ be as above. Then $M$ carries a
smooth foliation transverse to the Seifert fibration if and only 
if one of the following holds:
\begin{itemize}
\item
$-k+2\leq e_0\leq -2$
\item
$e_0=-1$ and $\Ga(M)$ is realizable
\item
$e_0=-k+1$ and $\Ga(-M)$ is realizable
\end{itemize}
\qed\end{thm}

The following theorem is part of a result~\cite[Theorem~1.3]{LM} of
the first named author and of Gordana Mati\'c (Ko Honda also obtained
similar results~\cite{H2}).

\begin{thm}[\cite{LM}]\label{t:quoteLM}
Let $M=M(e_0;r_1,r_2,\ldots,r_k)$ be as above. Then $M$ carries a positive 
contact structure transverse to the Seifert fibration if and only 
if one of the following holds:
\begin{itemize}
\item
$k\leq 2$ and $e(M):=e_0+\sum_i r_i <0$
\item
$e_0\leq -2$
\item
$e_0=-1$ and $\Ga(M)$ is realizable.
\end{itemize}
\qed\end{thm}

Recall that the Seifert fibered 3--manifold $M=M(e_0; r_1, r_2,\ldots,
r_k)$ is defined by the surgery diagram of Figure~\ref{f:seifert}.
Applying inverse slam--dunks, the diagram can be turned into a surgery
diagram involving only integer coefficients: replace each
$(-\frac{1}{r_i})$--surgery with a sequence of integer surgeries along
a chain of unknots, where the integral surgery coefficients are given
by the coefficients of the continued fraction expansion
\[
-\frac{1}{r_i} = a_0^{(i)} - \cfrac{1}{a_1^{(i)} - \cfrac{1}{\ddots -
    \cfrac{1}{a_k ^{(i)}} }},\quad a^{(i)}_j \leq -2.
\]
The integral surgery diagram defines a 4--manifold $W=W(e_0; r_1,
r_2,\ldots, r_k)$.  Notice that in the diagram all knots are unknots,
they are arranged along a star shaped tree $T$, and all framings are $\leq -2$
except the central one, which is equal to $e_0$. The 4--manifold $W$
is endowed with an $S^1$--action obtained by equivariant plumbing
according to $T$, where each vertex corresponds to a disk bundle over
a sphere. By~\cite[Chapter~2]{Or} $M$ is isomorphic (as a
$3$--manifold with $S^1$--action) to the boundary of $W$. Moreover, it
follows from the results of~\cite{McW} that the $S^1$--manifold $W$
carries a symplectic form $\om$ such that every orbit of the
$S^1$--action on $\del W=M$ is tangent to the kernel of $\om|_{\del
W}$. Recall that a {\sl symplectic filling} of a closed contact
3--manifold $(M,\xi)$ is a symplectic 4--manifold $(X,\om)$ such that
(i) $X$ is oriented by $\om\wedge\om$, (ii) $\del X=M$ as oriented
manifolds, and (iii) $\om|_\xi\neq 0$ at every point of $M$. An
immediate corollary of this discussion is the following:

\begin{prop}\label{p:filling}
$W=W(e_0; r_1, r_2,\ldots, r_k)$ carries a symplectic form $\om$ such
that if $M=M(e_0; r_1, r_2,\ldots, r_k)$ as above carries a positive
transverse contact structure $\xi$ then $(W,\om)$ is a symplectic
filling of $(M,\xi)$.\qed
\end{prop} 

It is known~\cite{NR} that
\[
b_2^+(W)=
\begin{cases}
1\quad\text{if}\quad e(M)>0,\\
0\quad\text{if}\quad e(M)<0.
\end{cases}
\]

The following result is not explicitely stated in~\cite{LM}, but it is
implicitely contained in the proof of~\cite[Theorem~1.3]{LM} given
in~\cite[Section~4]{LM}:

\begin{thm}[\cite{LM}]\label{t:plm}
Let $M=M(e_0;r_1,r_2,\ldots,r_k)$ be as above. If $M$ carries no positive
contact structure transverse to the Seifert fibration then there
exists an oriented surface $\Si$ smoothly embedded in $W$ such 
that $g(\Si)>0$ and 
\[
\Si\cdot\Si > 2g(\Si)-2.
\]
\qed
\end{thm}

\sh{Contact surgery}

Suppose that $L\subset (Y, \xi )$ is a Legendrian knot, that is, a
knot tangent to the 2--plane field $\xi $.  Oriented normals
to $\xi$ along $L$ provide a framing for $L$, called the~\emph{contact
framing}.  Let $Y_r(L)$ denote the result of $r$--surgery along $L$,
where the surgery coefficient is measured with respect to the contact
framing. Contact structures on $Y_r(L)$ can be defined by taking the
restriction $\xi \vert _{Y-\nu L}$ of $\xi $ to the complement of a
standard convex neighborhood of $L$, and extending it to a contact
structure on $Y_r(L)$ which is tight on the glued--up solid torus. By
the classification of tight contact structures on $S^1\times D^2$
\cite{Gi, H1}, such an extension exists if and only if $r\neq 0$ and is
unique if and only if $r=\frac{1}{k}$ for some $k\in \Z$. In general,
there are many tight structures on $S^1\times D^2$ satisfying the same
boundary condition. As Honda showed, a tight contact $(S^1\times D^2,
\eta )$ can be decomposed as
\[
(S^1\times D^2, \eta )=\cup _ {i=1}^t(T^2\times I, \eta _i ) \cup
(S^1\times D^2, \eta _0),
\]
where $(S^1\times D^2, \eta _0)$ is the standard neighborhood of a
Legendrian knot isotopic to $S^1\x 0\subset S^1\x D^2$ and each $(T^2\times
I, \eta _i)$ is a `continued fraction block'. A continued fraction
block $(T^2\times I, \eta _i)$ admits a further slicing
\[
(T^2\times I, \eta _i)= \cup _{j=1} ^{s_i} (T^2\times I, \eta _i (j)),
\]
where each $(T^2\times I, \eta _i (j))$ is a `basic slice' (see
\cite{H1} for definitions). Basic slices are
characterized by a sign $+$ or $-$. The continued fraction
block--decomposition of $(S^1\times D^2, \eta )$ is dictated by the
boundary value (the 'slope' of $T^2=\partial (S^1\times D^2)$), which,
in the case of surgery is determined by the surgery coefficient
$r$. The signs of the basic slices, however, rely on choices, giving
rise to possibly many tight contact structures on $S^1\times D^2$ with
identical boundary condition.

Ding and Geiges~\cite{DG2} showed how to realize each basic slice
decomposition of $(S^1\times D^2, \eta )$ by a contact surgery diagram
(cf. also \cite{DGS}). Let $r=\frac{p}{q}$ and if $r>0$ choose
$k\in\N$ minimal such that $r'=\frac{p}{q-kp}<0$. Consider $k$
Legendrian push--offs $L_1, \ldots , L_k$ of $L$ and perform contact
$(+1)$--surgeries along these push--offs. (In this part of the
procedure there is no choice.) Next, do contact $r'$--surgery along
$L$. Suppose that the contact solid torus $(S^1\times D^2, \eta )$ to
be used in the surgery has basic slice decomposition
\[
\cup _{i=1}^t \cup _{j=1} ^{s_i} (T^2\times I, \eta _i (j))\cup (S^1\times 
D^2, \eta _0).
\]
After fixing an orientation for $L$, apply a right stabilizations to
$L$ for each positive basic slice in
\[
(T^2\times I, \eta _1)=\cup _{j=1}^{s_1}(T^2\times I , \eta _1 (j))
\]
and a negative stabilization for each negative basic slice.  Denote
the result by $L(1)$. Consider a Legendrian push--off $L'(1)=L(2)$ of
$L(1)$ and perform contact $(-1)$--surgery along $L(1)$. Repeat the
above procedure for $L(2)$, using the second continued fraction block
of the decomposition. After $t$ steps the procedure terminates
providing a Legendrian link in $(Y, \xi )$ along which the
prescribed contact $(\pm 1)$--surgeries give the contact structure on
$Y_r(L)$ obtained by applying contact $r$--surgery with the prescribed
extension $(S^1\times D^2, \eta )$ on the glued--up solid torus.

\sh{Ozsv\'ath--Szab\'o homologies}

In the seminal papers \cite{OSzF1, OSzF2} a collection of homology
groups --- the Ozsv\'ath--Szab\'o homologies --- $\hf (Y, \t),
HF^{\pm}(Y,t)$ and $HF ^{\infty }(Y,\t)$ have been assigned to any
closed, oriented spin$^c$ 3--manifold $(Y, \t)$. For simplicity,
\emph{throughout this paper we shall use $\Z/2\Z$--coefficients}. The
groups in question admit a relative $\Z/d(\t)\Z$--grading, where
$d(\t)$ is the divisibility of the first Chern class $c_1(\t)$. When
$c_1(\t)$ is torsion, then  $d(\t)=0$, and the relative $\Z$--grading lifts
to an absolute $\Q$--grading.  The groups $HF^{\pm}(Y, \t)$ and
$HF^{\infty }(Y, \t)$ admit $\Z [U]$--module structures, and
multiplication by $U$ decreases grading by 2.

A spin$^c$ cobordism $(X, \s )$ from $(Y_1, \t _1)$ to $(Y_2, \t _2)$
induces a $\Z [U]$--equivariant homomorphism $F_{W, \s}$ between the
corresponding groups, and if $\t _1, \t _2$ are both torsion spin$^c$
structures then $F_{W, \s}$ shifts degree by
\[
\frac{1}{4}(c_1^2(\s ) -3\sigma (W)-2 \chi (W)).
\]
The direct sum of Ozsv\'ath--Szab\'o homology groups for all spin$^c$
structures is usually denoted by $\hf (Y) $ (respectively by
$HF^{\pm}(Y), \ HF^{\infty }(Y)$), while the sum of the maps induced by
$(X, \s )$ for all $\s \in Spin ^c (X)$ is denoted by $F_X$
(respectively $F^{\pm}_X,\  F^{\infty }_X$).

A rational homology sphere $Y$ is called an \emph{$L$--space} if $\hf
(Y, \t )$ is isomorphic to $\Z /2\Z $ for all spin$^c$ structure $\t
\in Spin ^c(Y)$. Equivalently, the dimension $\dim _{\Z /2\Z}\hf (Y) $ 
is equal to the order $\vert H_1(Y; \Z )\vert$.  Other equivalent ways
to define $L$--spaces is to require that the $U$--action on $HF^+(Y,
\t)$ is surjective, or that the natural map $HF^{\infty}(Y, \t )\to
HF^+(Y, \t)$ is onto.  Yet another characterization is that for an
$L$--space $Y$ the group $\hf (Y ,\t )$ can be identified with the
kernel of the $U$--map
\[
U\colon HF^+(Y, \t )\to HF^+(Y, \t).
\]
(The equivalences follow from the long exact sequences connecting the
various groups, see~\cite{OSzF2}.)  For a rational homology 3--sphere
$Y$ and spin$^c$ structure $\t \in Spin^c(Y)$ let the
\emph{$d$-invariant} $d(Y,\t)$ of $(Y,\t )$ be defined as the absolute
degree of the unique nontrivial element $x\in HF^{+}(Y,\t)$ which is
in the image of the natural map $HF^{\infty} (Y,\t)\to HF^{+}(Y, \t )$
and $Ux=0$. It is easy to see that for an $L$--space $Y$ the absolute
degree of the generator of $\hf (Y, \t)$ is equal to $d(Y, \t )$.  It
is known that $d(-Y, \t )=-d(Y, \t )$.

Suppose that the manifold $Y(K)$ is given as integral surgery along a
knot $K$ in $Y$, while $Y_1(K)$ is defined by an integral surgery
along $K\subset Y$ with framing one higher. According to 
\cite[Theorem~9.16]{OSzF2} these groups (together with the maps induced by
appropriate cobordisms $W_1, W_2, W_3$ between the 3--manifolds) fit into 
the exact triangle:
\[
\begin{graph}(6,2)
\graphlinecolour{1}\grapharrowtype{2}
\textnode {A}(1,1.5){$\hf (Y)$}
\textnode {B}(5, 1.5){$\hf (Y(K))$}
\textnode {C}(3, 0){$\hf (Y_1(K))$}
\diredge {A}{B}[\graphlinecolour{0}]
\diredge {B}{C}[\graphlinecolour{0}]
\diredge {C}{A}[\graphlinecolour{0}]
\freetext (3,1.8){$F_{W_1}$}
\freetext (4.6,0.6){$F_{W_2}$}
\freetext (1.4,0.6){$F_{W_3}$}
\end{graph}
\]

Ozsv\'ath--Szab\'o homologies are quite hard to compute in general,
but for 3--manifolds which can be presented as boundaries of plumbings
along negative definite plumbing trees with no `bad' vertices
(cf. \cite{OSzplum}), such computation has been carried out in
\cite{OSzplum}.  This immediately implies, for example, that a Seifert
fibered rational homology 3--sphere with $k$ multiple fibers and
$e_0(M)\leq -k$ is an $L$--space, since it can be presented as the
boundary of a plumbing tree without `bad' vertices. Since
$e_0(-M)=-k-e_0(M)$ for Seifert fibered spaces, the fact $\dim \hf
(Y)=\dim \hf (-Y)$ implies that \emph{a Seifert fibered 3--manifold
$M=M(e_0; r_1, r_2, \ldots , r_k)$ with $e_0(M)\leq -k$ or $e_0(M)\geq
0$ is an $L$--space}.

The following two lemmas can be deduced from the results
of~\cite{OSzplum}. They will be used in Section~\ref{s:firstproof}.

\begin{lem}\label{l:shift}
Suppose that $W$ is a negative definite cobordism between the
$L$--space $Y_1$ and the rational homology 3--sphere $Y_2$ obtained by
attaching 2--handles to $Y_1$. Let $\s\in Spin^c(W)$ be a spin$^c$
structure such that the degree shift of the maps
\[
F_{W, \s}\co\hf(Y_1,\s\vert_{Y_1})\to\hf (Y_2,\s\vert_{Y_2})
\]
and 
\[
F^+_{W,\s}\co HF^+(Y_1,\s\vert_{Y_1})\to HF^+(Y_2,\s\vert_{Y_2})
\]
is
\[
d(Y_2,\s\vert_{Y_2})-d(Y_1,\s\vert_{Y_1}).
\]  
Then the maps $F_{W, \s}$ and $F^+_{W,\s}$ are injective.
\end{lem}

\begin{proof}
Since $W$ is negative definite, successive applications
of~\cite[Proposition~9.4]{abs} imply that the map
\[
F_{W,\s}^{\infty}\colon HF^{\infty }(Y_1, \s\vert _{Y_1})\to 
HF^{\infty }(Y_2, \s \vert _{Y_2})
\]
is an isomorphism for any spin$^c$ structure. This implies that the
map $F^+_{W, \s}$ is injective provided its degree shift is equal to
the difference of the $d$-invariants. The long exact sequence
connecting $\hf (Y, \t )$ and $HF ^+(Y, \t )$ together with the fact
that the natural map
\[
\hf(Y_1, \s \vert _{Y_1})\to HF^+(Y_1,\s \vert _{Y_1})
\]
is injective for the $L$--space $Y_1$ imply that if $F^+_{W,\s}$ is
injective then so is the map $F_{W, \s }$ between the $\hf$--groups,
concluding the proof.
\end{proof}

\begin{lem}\label{l:spincs}
Suppose that $W$ is the cobordism from $S^3$ to an $L$--space $L$
given by a negative definite plumbing with no `bad' vertex in the
sense of~\cite{OSzplum}. Then for any spin$^c$ structure $\t \in
Spin^c (L)$ there is a spin$^c$ structure $\s \in Spin^c(W)$ such that
the map
\[
F_{W,\s}\colon \hf (S^3) \to \hf (L, \t)
\]
is an isomorphism.
\end{lem}

\begin{proof}
The proof is straightforward using the algorithm of~\cite{OSzplum}
which computes the kernel of the $U$--map on $HF^+(L, \t )$. Take the
negative definite plumbing $W$ for $L$ and choose the characteristic
vector $K\in H^2 (W; \Z )$ with the following three properties:
\begin{itemize}
\item $K$ satisfies the starting condition given in the
first paragraph of~\cite[Section~3.1]{OSzplum},
\item $K$ leads to a final vector $L$ with the property given in
$(16)$ of the same paragraph, and 
\item $K$ induces the spin$^c$ structure $\t$ on $L$.
\end{itemize}
Since the algorithm computes the Ozsv\'ath--Szab\'o homology of the
boundary $L$--space $L$ with the given spin$^c$ structure, such a
vector clearly exists. By the formula for the grading given
in~\cite{OSzplum}, the final vector of this process defines a spin$^c$
structure on the 4--manifold $W$ satisfying the conditions of
Lemma~\ref{l:shift}. Since $L$ is an $L$--space, this clearly suffices
to prove the statement.
\end{proof}

\sh{Contact Ozsv\'ath--Szab\'o invariants}

A contact structure $\xi $ on $Y$ determines an element $c(Y,\xi )\in
\hf (-Y, \t _{\xi })$ which has the following crucial
properties~\cite{OSzcont}:
\begin{itemize}
\item $c(Y, \xi )$ is an isotopy invariant of the contact 3--manifold
$(Y, \xi )$;
\item $c(Y, \xi )=0$ if the contact structure $\xi $ is overtwisted;
\item $c(Y, \xi )\neq 0$ if $(Y, \xi )$ is Stein fillable;
\item if $c_1(\t _{\xi })$ is a torsion class and $c(Y, \xi )\neq 0$
 then it is of degree $-d_3(\xi )$ in $\hf (-Y, \t _{\xi })$, where
 $d_3(\xi )\in\Q$ is the Hopf invariant of the 2--plane field underlying
 $\xi$.
\end{itemize}

\section{The proof of Theorem~\ref{t:equiv}}\label{s:firstproof}

Some of the implications among the equivalent statements of
Theorem~\ref{t:equiv} are relatively easy to prove (and have been
partly established in the literature), so we will start with those.

$(1)\Longrightarrow (2)$ Suppose that $M$ carries transverse contact
structures with both its orientations. Since $e(-M)=-e(M)$, we can
choose the orientation for which the 4--manifold $W$ given by the
integer plumbing representation of $M$ has positive
$b_2^+$--invariant.  According to Proposition~\ref{p:filling} this
manifold carries a symplectic structure which provides a symplectic
filling for any transverse contact structure. Since $b_2^+(W)>0$, by
\cite[Theorem~1.4]{OSzgen} $M$ cannot be an $L$--space.

$(2)\Longrightarrow (3)$ Since $M\neq S^1\x S^2$, according to
\cite{ET} a transverse foliation can be $C^0$--approximated by
positive as well as negative transverse contact structures. Therefore,
the nonexistence of positive transverse contact structures with one
orientation prevents the existence of a transverse foliation on $M$.

$(1)\Longrightarrow (4)$ This is part
of~\cite[Theorem~1.4]{OSzgen}.

$(4)\Longrightarrow (3)$
Since transverse foliations are taut, this implication is trivial.

$(3)\Longrightarrow (2)$ Suppose that $M=M(e_0;r_1,r_2,\ldots,r_k)$
admits no transverse foliation. Either $e(M)>0$ or $e(-M)=-e(M)>0$.
Therefore, if $k\leq 2$ the conclusion follows from
Theorem~\ref{t:quoteLM}. Now suppose that $k\geq 3$. According to
Theorem~\ref{t:quoteJN}, up to changing the orientation of $M$ we may
assume that either $e_0(M)\geq 0$ or $e_0(M)=-1$ and $\Ga(M)$ is not
realizable. In both cases, by Theorem~\ref{t:quoteLM} $M$ carries no
positive, transverse contact structures.

In the light of our observations above, to complete the proof of the
equivalences we only need to show that $(2)$ implies $(1)$, that is,
if the Seifert fibered rational homology 3--sphere $M$ carries no
transverse contact structures with one of its orientations then it is
an $L$--space. Observe that if $k\leq 2$ then $M$ is a lens space,
therefore in this case there is nothing to prove. Likewise, if either
$e_0(M)\leq -k$ or $e_0(M)\geq 0$, then $M$ is an
$L$--space. Moreover, if $-k+1\leq e_0(M)\leq -2$ then by
Theorem~\ref{t:quoteLM} $M$ carries transverse contact
structures. Therefore, up to changing the orientation of $M$ we may
assume that $e_0(M)=-1$, so to establish Theorem~\ref{t:equiv} it will
suffice to prove the following:

\begin{prop}\label{p:lastcase}
If the rational homology 3--sphere $M=M(-1;r_1,\ldots,r_k)$ 
carries no transverse, positive contact structures, then $M$ is an
$L$--space.
\end{prop}

Before proving Proposition~\ref{p:lastcase} we need some
preliminaries. Suppose that $M$ does not carry a transverse positive
contact structure.  By Theorem~\ref{t:plm}, this implies that the
associated 4--manifold $W$ contains a surface $\Sigma$ with
$g(\Sigma)>0$ and $\Si\cdot\Si > 2g(\Sigma)-2$. Let $W_0$ be
the 4--manifold obtained from $W$ by deleting an open 4--ball from its
interior. Then, by~\cite[Proposition~2.1]{LSgt} the map
\[
F_{W_0,\s}\co \hf (S^3)\lra \hf (M)
\]
vanishes for every spin$^c$ structure $\s\in Spin ^c(W_0)$.  Now let
us consider the subcobordism $W_1$ obtained by attaching to $S^3$ the
2--handles corresponding to the vertices of the plumbing tree with
framings at most $-2$. This 4--manifold is a cobordism between $S^3$
and a connected sum of $k$ lens spaces $L$. By standard properties
of lens spaces and connected sums, $L$ is an $L$--space. Therefore,
$W$ satisfies the assumptions of Lemma~\ref{l:spincs}. The attachment
of the $(-1)$--framed circle induces a cobordism $W_2$ from $-L$ to
$M$.  Since $W_1$ is negative definite and $W_0=W_1\cup W_2$ is not
(because it contains a surface with positive self--intersection), we
conclude that $b_2^+(W_2)=1$.

\begin{lem}\label{l:zero}
For each spin$^c$ structure $\s_2 \in Spin^c(W_2)$ the map $F_{W_2, \s
_2}$ on $\hf (L)$ is zero.
\end{lem}

\begin{proof}
Let $x\in \hf (L,\t )$ be a given generator and consider $\s _1\in
Spin ^c(W_1)$ provided by Lemma~\ref{l:spincs} with the property
that $F_{W_1, \s _1}(g)=x$ for the generator $g\in \hf (S^3)$.  By
the composition law~\cite[Theorem~3.4]{OSzF4} we have 
\[
F_{W_2, \s _2 }(x)=(F_{W_2,\s _2}\circ F_{W_1, \s _1})(g)= (\sum _{\t \in
{\mathfrak{S}}} F_{W_0, \t})(g)=0,
\]
where ${\mathfrak {S}}=\{ \t \in Spin ^c(W_0) \mid \t \vert _{W_i}=\s
_i, \ i=1,2\}$. Since $F_{W_0, \t}=0$ for all $\t \in Spin ^c (W_0)$,
this implies $F_{W_2, \s _2}=0$.
\end{proof}

\begin{proof}[Proof of Proposition~\ref{p:lastcase}]
To verify that $M$ is an $L$--space, consider the surgery exact
triangle induced by the cobordism $W_2$:
\[
\begin{graph}(6,2)
\graphlinecolour{1}\grapharrowtype{2}
\textnode {A}(1,1.5){$\hf (-L)$}
\textnode {B}(5, 1.5){$\hf (M)$}
\textnode {C}(3, 0){$\hf (M')$}
\diredge {A}{B}[\graphlinecolour{0}]
\diredge {B}{C}[\graphlinecolour{0}]
\diredge {C}{A}[\graphlinecolour{0}]
\freetext (3,1.8){$F_{W_2}$}
\end{graph}
\]
Notice that $M'$ is a small Seifert fibered space with $e_0(M')=0$,
hence $M'$ is an $L$--space. Since $F_{W_2}=0$, we have
\[
\hf (M')=\hf (L)\oplus \hf (M).
\]
According to~\cite[Theorem~8.1]{NR}, if $Y=M(e_0;r_1,r_2,\ldots,r_k)$
with with $r_i=\frac{p_i}{q_i}$, $i=1,2,\ldots,k$ and $e(Y)\neq 0$,
then the order of $H_1(Y;\Z)$ is equal to
\[
q_1 q_2\cdots q_k |e(Y)|. 
\]
Therefore, since $e(M)>0$,  
\[
\begin{split}
|H_1(M;\Z)| & = (\prod_{i=1}^k q_i)(\sum_{i=1}^k \frac{p_i}{q_i}-1) 
= \sum_{i=1}^k p_i\prod_{j\neq i} q_j - \prod_{i=1}^k q_i \\
& = |H_1(M';\Z)| - |H_1(L;\Z)|.
\end{split}
\]
Thus, we obtain 
\[
\dim \hf (M)=\vert H_1(M; \Z )\vert, 
\]
concluding the proof.
\end{proof}

\section{The proof of Theorem~\ref{t:suff}}\label{s:tight}

\begin{figure}[ht]
\begin{center}
\psfrag{L1}{\small $L_1$}
\psfrag{L2}{\small $L_2$}
\psfrag{Lk}{\small $L_k$}
\psfrag{K}{\small $K_0$}
\psfrag{+1}{\small $+1$}
\psfrag{r1}{\small $-\frac1{r_1}$}
\psfrag{r2}{\small $-\frac1{r_2}$}
\psfrag{rk}{\small $-\frac1{r_k}$}
\epsfig{file=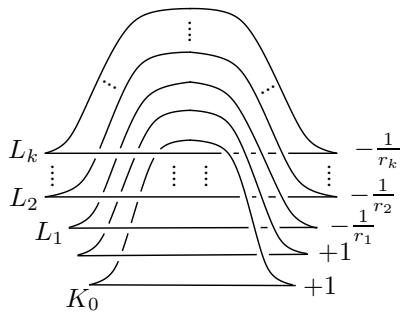, height=4cm}
\end{center}
\caption{Contact structures on $M$}
\label{f:structures}
\end{figure}

Let $M=M(-1; r_1, r_2,\ldots, r_k)$ be a 3--manifold as in the statement of
Theorem~\ref{t:suff}, satisfying $e(M)=r_1+r_2+\cdots + r_k-1>0$. Define the
family of contact 3--manifolds $\{(S, \xi_S)\}$ using the contact
surgery diagram of Figure~\ref{f:structures} with the knot $K_0$ deleted.

\begin{lem}\label{l:S}
Each contact 3--manifold $(S,\xi_S)$ is Stein fillable, and each underlying 
3--manifold $S$ is an $L$--space.
\end{lem}

\begin{proof}
Let $K$ be a Legendrian unknot with maximal Thurston--Bennequin number
in the standard contact 3--sphere. A contact $(+1)$--surgery along $K$
gives a Stein fillable contact 3--manifold $(S^1\x S^2,\eta)$. This
follows e.g.~from the fact that the resulting contact structure has
nonvanishing Ozsv\'ath--Szab\'o invariant~\cite[Lemma 2.5]{LS1.5} and
the classification of tight contact structures on $S^1\x S^2$.
Since $(S, \xi _S)$ is given by Legendrian surgery on
$(S^1\x S^2,\eta)$, $\xi _S$ is Stein fillable.  Since $e_0(M)=-1$, a
simple calculation shows that $S$ is a small Seifert fibered
3--manifold with $e_0(S)=0$, therefore it is an $L$--space.
\end{proof}

Let the 4--manifold $-X$ be defined as the cobordism induced by the
contact $(+1)$--surgery along the Legendrian unknot $K_0$ in
Figure~\ref{f:structures}. Note that $-X$ is a cobordism between $S$
and $M$. By reversing its orientation, the resulting 4--manifold (a
cobordism between $-S$ and $-M$) will be denoted by $X$. 

\begin{lem}\label{l:cobo}
The cobordism $X$ is negative definite.
\end{lem}

\begin{proof}
Converting the contact framings into smooth framings in
Figure~\ref{f:structures}, reversing orientation, blowing up once and
applying three Rolfsen twists one easily checks that the cobordism $X$
is given by Figure~\ref{f:W}, where the surgery presentation for $-S$ is
drawn using solid curves, and the attaching circle of the 2--handle
giving $X$ is drawn using the dashed circle.
\begin{figure}[ht]
\begin{center}
\psfrag{r1}{\small $-\frac1{1-r_1}$}
\psfrag{r2}{\small $-\frac1{1-r_2}$}
\psfrag{rk}{\small $-\frac1{1-r_k}$}
\psfrag{-1}{\small $-1$}
\psfrag{-3}{\small $-k$}
\epsfig{file=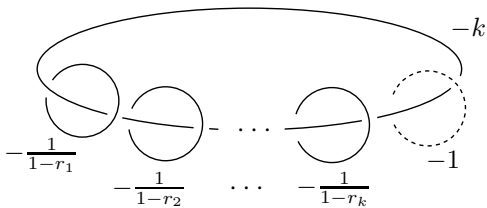, height=2.5cm}
\end{center}
\caption{The cobordism $X$}
\label{f:W}
\end{figure}
Using a little bit of Kirby calculus, one can easily see that the
cobordism $X$ admits an embedding into the 4--manifold
\[
\widehat W = W(-k+1; 1-r_1,1-r_2,\ldots,1-r_k)\# \cpkk. 
\]
Since 
\[
\del\widehat W = M(-k+1;1-r_1,1-r_2,\ldots,1-r_k) = -M
\]
and $e(-M) = -e(M) < 0$, we conclude that $W(-k+1; 1-r_1,1-r_2,\ldots,1-r_k)$
is negative definite, hence so is $\widehat W$, consequently the
cobordism $X$ is negative definite.
\end{proof}

\begin{proof}[Proof of Theorem~\ref{t:suff}]
Denote by $c^+(M, \xi)$ the image of $c(M, \xi)$ under the natural
homomorphism~\cite{OSzF1, OSzF2}
\[
\hf(-M,\t_{\xi})\to HF^+(-M,\t_{\xi}).
\]
Clearly, it is enough to show that $c^+(M, \xi)\neq 0$. The contact
3--manifold $(M,\xi)$ is obtained by a contact $(+1)$--surgery on
$(S,\xi_S)$, for some $(S,\xi_S)$ from the family of
Lemma~\ref{l:S}. Therefore, by~\cite[Lemma~2.11]{Gh}, there is a
spin$^c$ structure $\s$ on $X$ such that
$F^+_{X,\s}(c^+(S, \xi _S))=c^+(M,\xi)$ and
\[
-d_3(\xi_S) + \de(\s) = -d_3(\xi),
\]
where
\[
\de(\s):=\frac 14 (c_1^2(\s) - 3\si(-X) - 2\chi(-X)).
\]
By assumption, $M$ is a rational homology 3--sphere, and by
Lemma~\ref{l:S} $S$ is an $L$--space and $(S, \xi_S)$ is Stein
fillable. Therefore $c(S, \xi_S)\neq 0$, and it follows that
$d(-S,\t_{\xi_S})=-d_3(\xi_S)$. By Lemma~\ref{l:cobo} $X$ is negative
definite. Moreover, the assumption $d_3(\xi)=d(M, \t_{\xi})$
immediately implies that the degree shift of the map $F^+_{X,\s}$ is
\[
\de(\s)=-d_3(\xi) + d_3(\xi_S) = d(-M,\t_{\xi}) - d(-S,\t_{\xi_S}).
\]
Therefore, Lemma~\ref{l:shift} applies and the map $F^+_{X,\s}$ is
injective, so we conclude that $c^+(M, \xi)\neq 0$.
\end{proof}

\section{The proof of Theorem~\ref{t:existence}}\label{s:critical}

Let $K\subset S^3$ be a knot, and denote by $S^3_r(K)$ the 3--manifold
obtained by $r$--surgery along $K$. The proof of
Theorem~\ref{t:existence} will rest on Theorem~\ref{t:suff}, together
with Propositions~\ref{p:spin}, \ref{p:compd}, \ref{p:alex} and
Lemma~\ref{l:d-inv}. In the following $M_p$ will denote the
3--manifold $S^3_{p^2-p-1}(T_{p,p+1})$.

\begin{prop}\label{p:spin}
Let $\t_0$ denote the spin$^c$ structure induced by the (unique) spin
structure on the $L$--space $M_p$. Then 
\[
d(M_p, \t_0) =
\begin{cases}
-\frac{1}{4}(3p+2)\quad\text{if $p$ is even},\\
-\frac{1}{4}(p+1)\quad\text{if $p$ is odd}.
\end{cases}
\]
\end{prop}

\begin{proof}
If $p$ is even then $-M_p$ is obtained by plumbing according to a
negative definite tree with one bad vertex in the sense
of~\cite{OSzplum}, and with $3p+2$ vertices having all even
weights (see e.g.~\cite[Figure~4]{LS2}, where $M_p=E_{p,1}$). 

Therefore, the trivial vector $K=0$ is characteristic and clearly
induces $\t_0$, so by~\cite[Corollary~1.5]{OSzplum} we have
\[
-d(M_p,\t_0)= d(-M_p,\t_0)=\frac{K^2+3p+2}{4} = \frac{3p+2}4.
\]
If $p$ is odd, the proof follows easily from~\cite{LS2}.  Namely,
since $M_p=E_{p,1}$ and $\t_0=\t_E$ in the notation of~\cite{LS2},
by~\cite[Lemma~6.3]{LS2} and~\cite[Proposition~6.13]{LS2} we have 
\[
d(-M_p,\t_0)=d(-L_p,\t_Z)+\frac14,
\]
where $-L_p$ is the oriented boundary of the smooth 4--manifold
$Z_p$ given by attaching 2--handles to $B^4$ along the framed link
of Figure~\ref{f:Zp}, 
\begin{figure}[h]
\begin{center}
\psfrag{-2}{\small $-2$}
\psfrag{p-1}{\small $p-1$}
\psfrag{2p}{\small $2p$}
\psfrag{2}{\small $2$}
\psfrag{p+1}{\small $p+1$}
\epsfig{file=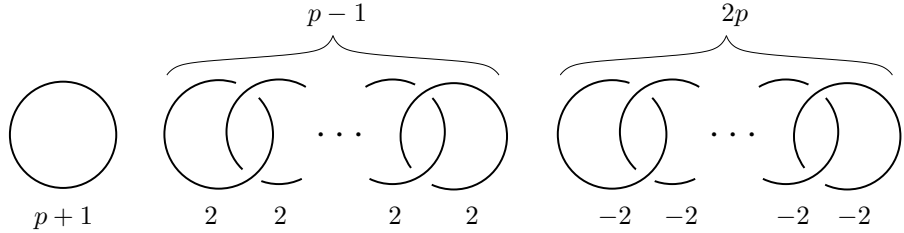, width=0.9\textwidth}
\end{center}
\caption{The framed link defining $Z_p$}
\label{f:Zp}
\end{figure}
and $\t_Z$ is the spin$^c$ structure induced by the restriction of the
unique spin structure on $Z_p$. A simple computation gives
\[
d(-L_p, \t_Z) = \frac{p}{4},
\]
and the conclusion follows.
\end{proof}

Let $X_p$ denote the 4--manifold obtained by attaching a 2--handle to
$D^4$ along $T_{p,p+1}\subset S^3$ with framing $p^2-p-1$. In
particular, $\partial X_p=M_p$. Choose a generator 
\[
[\Sigma_p]\in H_2(X_p; \Z )\cong\Z. 
\]
Since $X_p$ is simply connected, a spin$^c$ structure on $X_p$ is
uniquely determined by its first Chern class, so let 
$\s_k\in Spin^c(X_p)$ denote the spin$^c$ structure with
\[
\langle c_1(\s _k), [\Sigma _p ]\rangle = p^2-p-1+2k.
\]
We denote the restriction of $\s _k$ to $\partial X_p=M_p$ by $\t
_k$. Observe that, in accordance with the notation used in
Proposition~\ref{p:spin}, $\t_0$ is the spin$^c$ structure induced by the
unique spin structure on $M_p$. Let
\[
\Delta _{T_{p,p+1}}(t)=a_0(p)+\sum _{i=1}^{n_p} a_i (p) (t^i +t^{-i})
\]
be the symmetrized Alexander polynomial of $T_{p,p+1}$, and define
its $j^{th}$ torsion coefficient by
\[
r_j (p)=\sum _{i=1}^{n_p} ia_{i+\vert j \vert }(p).
\]

\begin{prop}\label{p:compd} 
For $|k|\leq\frac{p^2-p-1}2$, the invariant $d(M_p, \t_k)$ is equal
to
\[ 
\frac{((p^2-p-1)-2j )^2-(p^2-p-1)}{4(p^2-p-1)} - 2r_k(p),  
\]
where $j\equiv k\bmod (p^2-p-1)$ and $0\leq j<p^2-p-1$. 
\end{prop}

\begin{proof} 
Recall that $S^3_{pq\pm 1}(T_{p,q})$ are lens spaces. Therefore,
combining \cite[Theorem~7.2]{abs} (together with
\cite[Remark~7.4]{abs}) and \cite[Proposition~8.1]{abs}
(cf.~also~\cite[Theorem~1.2]{OSrat}) we have
\[
d(M_p, \t_k) - d(L(p^2-p-1,1),k) = -2r_k(p).
\]
On the other hand, by~\cite[Section~7]{abs} 
\[
d(L(p^2-p-1,1),k) = \frac{((p^2-p-1)-2j )^2-(p^2-p-1)}{4(p^2-p-1)},
\]
where $j\equiv k\bmod (p^2-p-1)$ and $0\leq j<p^2-p-1$. The statement 
follows immediately.
\end{proof}

\begin{prop}\label{p:alex} 
The coefficients of the Alexander polynomial $\Delta
_{T_{p,p+1}}(t)$ satisfy:
\begin{enumerate}
\item[(i)]
$a_0(p)=(-1)^{p+1}$,
\item[(ii)]
all the $a_i(p)$'s are $\pm 1$ and they alternate in sign, and
\item[(iii)]
$a_i(p)=0$ for $0<i <\frac{p}{2} $.
\end{enumerate}
\end{prop}

\begin{proof}
It is known (see for example~\cite[Example~9.15]{Bu}) that
\[ 
\Delta _{T_{p,p+1}}(t)=
\frac{(1-t^{p(p+1)})(1-t)}{(1-t^p)(1-t^{p+1})}t^{-\frac{p(p-1)}{2}}.\]
Since
\[ 
(1-t^{p(p+1)})=(1-t^p)(1+t^p+t^{2p}+ \dots +t^{p^2}),
\]
and
\[
\begin{split} 
& (1+t^p+t^{2p}+\dots+t^{p^2})(1-t) = \\
= & (1-t^{p+1})+(t^p-t^{2p+1})+\dots+(t^{p^2-p}-t^{p^2+1})+t^{p^2}-t=\\
= & (1-t^{p+1})(1+t^p+\dots+t^{p^2-p})+t(t^{p^2-1}-1)=\\
= & (1-t^{p+1})(1+t^p+\dots+t^{p^2-p})-\\
& t(1-t^{p+1})(1+t^{p+1}+t^{2(p+1)}+\dots+t^{(p+1)(p-2)}),
\end{split}
\]
we get that the (symmetrized) Alexander polynomial $\Delta
_{T_{p,p+1}}$ equals
\[
q(t)t^{-\frac{p(p-1)}{2}},
\]
where
\begin{equation}\label{e:q}
q(t) =1+t^p+\dots+t^{p^2-p}-t(1+t^{p+1}+t^{2(p+1)}+
\dots+t^{(p+1)(p-2}).
\end{equation}
Therefore, $a_0(p)$ is the coefficient of $t^{\frac{p(p-1)}{2}}$ in
$q(t)$, which is equal to $1$ if $\frac{1}{2}p(p-1)$ is of the form
$i\cdot p$ with $1\leq i \leq p-1$ (this holds if $p$ is odd), and it
is equal to $-1$ if $\frac{1}{2}p(p-1)$ is of the form $1+k(p+1)$
(which happens with the choice of $k=\frac{1}{2}(p-2)$ for $p$
even). This verifies $(i)$. Part $(ii)$ follows immediately from a
close inspection of Formula~\eqref{e:q}. 
To verify $(iii)$, suppose first that $p$ is
odd. The first coefficient after $a_0(p)$ which is equal to 1 is
clearly $a_p(p)$, while the first $-1$ is of index $1+k$ with
$k=\frac{1}{2}(p-1)$, showing that $a_{\frac{1}{2}(p-1) +1}(p)=
a_{\frac{1}{2}(p+1)}(p)$ is the first nonzero coefficient after
$a_0(p)$, verifying $(iii)$ for odd $p$. For $p$ even, a similar
argument shows that the first coefficient which is equal to $-1$ is
$a_{p+1}(p)$, and the first $1$ is $a_{\frac{p}{2}}(p)$, completing
the proof.
\end{proof}

\begin{lem}\label{l:d-inv}
For $p>2$ even, we have 
\[ 
d(M_p, \t _{\frac{p}{2}})= \frac{(p^2-2p-1)^2}{4(p^2-p-1)}-\frac{1}{4}(p-1)^2.
\] 
\end{lem} 

\begin{proof}
According to Proposition~\ref{p:spin} we have 
\[
d(M_p, \t_0)=-\frac{3p+2}{4},
\] 
and by Proposition~\ref{p:compd}   
\[
d(M_p, \t_0)=\frac{p^2-p-2}{4}-2r_0(p).  
\] 
The above identities imply that 
\[ 
r_0(p)=\frac{3p+2}{8}+\frac{p^2-p-2}{8}=\frac{p^2+2p}8.  
\] 
On the other hand, by Proposition~\ref{p:spin}
\[ 
d(M_p,\t_\frac{p}{2})
=\frac{(p-(p^2-p-1))^2-(p^2-p-1)}{4(p^2-p-1)}-2r_{\frac{p}{2}}(p).
\] 
To conclude the proof we will express $r_{\frac{p}{2}}(p)$ in 
terms of $r_0(p)$.  
It follows from the definition that 
\[
r_0(p)=r_{\frac{p}{2}}(p)+\frac{p}{2}\sum _{i=0}^{n_p} a_i(p)-\sum
_{i=0}^{\frac{p}{2}}(\frac{p}{2}-i)a_i(p).
\]
Since by Proposition~\ref{p:alex} the nonzero coefficients of $\Delta
_{T_{p,p+1}}(t)$ alternate in sign and are all $\pm 1$, from
the assumption that $p$ is even (hence $a_0(p)=-1$) we get 
\[
\sum _{i=0} ^{n_p} a_i(p)=0.
\]
On the other hand, by Proposition~\ref{p:alex} we know that 
\[
\sum _{i=0}^{\frac{p}{2}}(\frac{p}{2}-i)a _i (p)=-\frac{p}{2},
\]
hence
\[
r_{\frac{p}{2}}(p)=r_0(p)-\frac{p}{2}= \frac{p^2-2p}{8}.
\]
Substituting this value into the expression for $d(M_p, \t
_{\frac{p}{2}})$ given above, the statement follows.
\end{proof}

\begin{proof}[Proof of Theorem~\ref{t:existence}]
Suppose first that $p$ is even. Define $\xi _p$ by the contact surgery
diagram of Figure~\ref{f:pp+1even}. (In Figures~\ref{f:pp+1even}
and~\ref{f:pp+1odd} a coefficient $(+1)$ next to a Legendrian knot $K$
means that contact $(+1)$--surgery is performed along $K$, while no
coefficient means contact $(-1)$--surgery).
\begin{figure}[ht]
\begin{center}
\psfrag{p/2stabilizationssssssssssssss}{\small $\frac p2$ stabilizations}
\psfrag{2p}{\small $2p$ stabilizations}
\psfrag{+1}{\small $+1$}
\psfrag{p-1}{\small $p-1$ unknots}
\epsfig{file=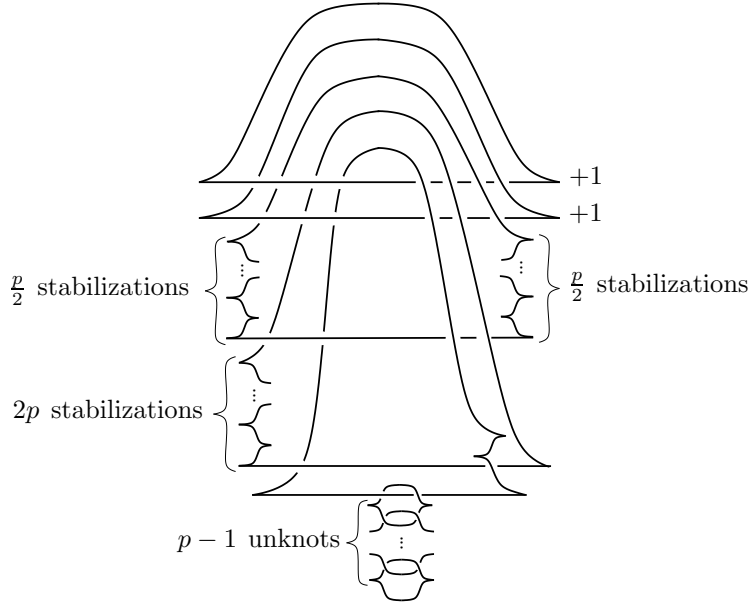, height=8cm}
\end{center}
\caption{The contact structure $\xi_p$ for $p$ even}
\label{f:pp+1even}
\end{figure}
Our plan is to apply Theorem~\ref{t:suff}. Observe that this plan
makes sense. In fact, it is easy to check that the proofs of
Lemma~\ref{l:S}, Lemma~\ref{l:cobo} and Theorem~\ref{t:suff} given in
Section~\ref{s:tight} apply without modifications to the contact
surgery presentation of Figure~\ref{f:pp+1even}.

Turning the contact framings into smooth, a little Kirby
calculus (as in~\cite[Figures~1 and~7]{LS2}) shows that $\xi _p$ is a
contact structure on $M_p$. Now we want to apply the formula
from~\cite{DGS} for the $d_3$--invariant of a contact structure
defined by a contact $(\pm 1)$--surgery diagram. If $c\in H^2(X;\Z)$ 
denotes the 2--cohomology class determined by the rotations numbers
(see~\cite{DGS}), $\si(X)$ is the signature of $X$ and $b_2(X)$ 
the second Betti number, a simple computation yields
\[
\sigma (X)=1-p,\quad b_2(X)=p+3,\quad 
c^2=\frac{(p-(p^2-p-1))^2}{p^2-p-1}-p(p-1).
\]
Then, \cite[Corollary~3.6]{DGS} (where $b_2(X)$ should be plugged into
the formula instead of the Euler characteristic $\chi(X)$ because the
3--dimensional invariant used in Heegaard Floer theory is shifted by
$\frac12$) gives
\[
d_3 (\xi _p) = \frac{(p^2-2p-1)^2}{4(p^2-p-1)}-\frac{1}{4}(p-1)^2.
\]
Following the blow--down procedure at the cohomological level, the
verification that $\xi _p$ is a contact structure on the given
3--manifold $M_p$ also shows that for $p>2$
\[
\t_{\xi _p}=\t _{\frac{p}{2}}.
\]
(For $p=2$ the inequality $p<p^2-p-1$ fails to hold, and we have
$\t_{\xi_2}=\t_0$.)

Therefore, by Lemma~\ref{l:d-inv} we have 
\[
d_3(\xi _p) = d(M_p,\t_{\xi _p}).
\]
Using Theorem~\ref{t:suff} we conclude that $c(M_p,\xi_p)\neq 0$,
hence $\xi_p$ is a tight contact structure on $M_p$.

We now verify the statement for $p$ odd. (This case was already
treated in~\cite[Theorem~1.3]{LS2}.)  Let $\xi_p$ denote the contact
structure given by the contact surgery diagram of
Figure~\ref{f:pp+1odd}.
\begin{figure}[ht]
\begin{center}
\psfrag{(p+1)/2}{\small $\frac {p+1}2$ stabilizations}
\psfrag{(p-1)/2stabilizationssssssssssssss}
{\small $\frac {p-1}2$ stabilizations}
\psfrag{2p}{\small $2p$ stabilizations}
\psfrag{+1}{\small $+1$}
\psfrag{p-1}{\small $p-1$ unknots}
\epsfig{file=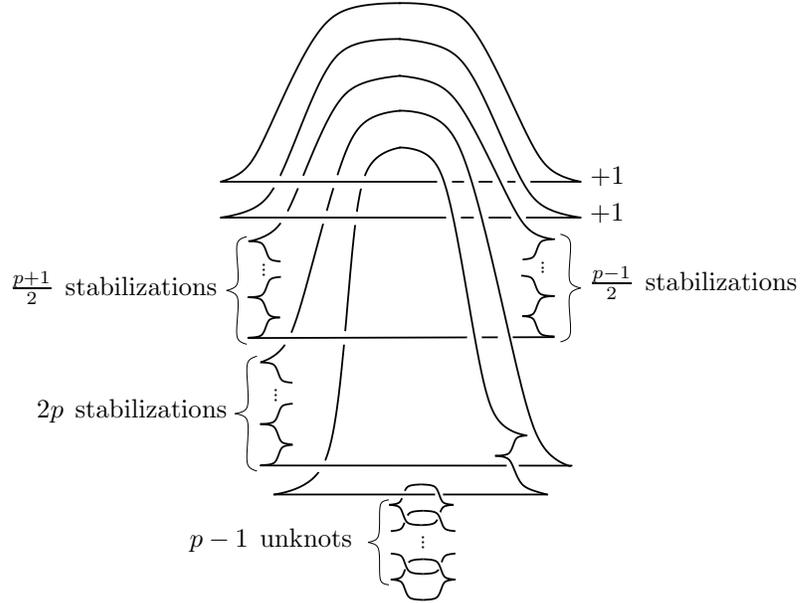, height=8cm}
\end{center}
\caption{The contact structure $\xi _p$ for $p$ odd}
\label{f:pp+1odd}
\end{figure}
As in the previous case, we can apply Theorem~\ref{t:suff}. A
computation as above shows that
\[
\sigma(X) = 1-p,\quad b_2(X)=p+3,\quad c^2=-2p
\]
and applying~\cite[Corollary~3.6]{DGS} we have
\[
d_3(\xi _p)=-\frac{1}{4}(p+1)\quad\text{and}\quad \t_{\xi_p}=\t _0.
\] 
Using Proposition~\ref{p:spin} we see that 
\[
d(M_p,\t_{\xi_p})=-\frac{1}{4}(p+1) = d_3(\xi _p)
\]
which, in view of Theorem~\ref{t:suff}, implies that $\xi_p$ is tight
and concludes the proof.
\end{proof}

\section{Planar contact structures}
\label{s:six}

\sh{Surgery diagrams for contact structures with zero twisting}

In this section we show that if $\xi$ is a tight contact structure
with maximal twisting equal to zero on the small Seifert fibered
3--manifold $M=M(-1; r_1, r_2, r_3)$, then $(M,\xi)$ is compatible
with a planar open book decomposition of $M$, and we discuss some
consequences. We start by establishing the fact that $(M,\xi)$ is
given by one of the contact surgeries represented by
Figure~\ref{f:structures}.

Let $F_i$ ($i=1,2,3$) be the three singular fibers of the Seifert
fibration on $M$. We first isotope $F_i$ so that they become
Legendrian.  Then, since $\xi$ has maximal twisting equal to zero, we
can isotope the Seifert fibration further, so that there is a
Legendrian regular fiber $L$ with contact framing equal to the framing
induced on $L$ by the Seifert fibration.

Let $V_i$ be a standard neighbourhood of $F_i$, $i=1,2,3$. Then $M
\setminus (V_1 \cup V_2 \cup V_3)$ can be identified with $\Sigma
\times S^1$ where $\Sigma$ is a pair--of--pants. An identification
between $M \setminus (V_1 \cup V_2 \cup V_3)$ and $\Sigma \times S^1$
determines identifications of $- \partial (M \setminus V_i)$ with
$\R^2 / \Z^2$ so that $\binom{1}{0}$ is the direction of the section
$\Sigma \times \{ 1 \}$ and $\binom{0}{1}$ is the direction of the
regular fibers. A standard application of convex surface theory (as in
e.g.~\cite{GLS2}) shows that the singular fibers $F_1, F_2, F_3$ admit
tubular neighbourhoods $U_1, U_2, U_3$ with minimal convex boundaries
such that $-\partial (M \setminus U_i)$ has infinite slope for
$i=1,2,3$. 

\begin{prop}\label{p:model}
Let $\xi$ be a tight contact structure with maximal twisting equal to
zero on the small Seifert fibered 3--manifold $M=M(-1; r_1, r_2,
r_3)$. Then $(M,\xi)$ is given by one of the contact surgeries
represented by Figure~\ref{f:structures}.
\end{prop}

\begin{proof}
Let $\Si$ be a pair of pants. Recall that a tight contact structure
$\xi$ on $\Sigma \times S^1$ is {\em appropriate} if there is no
contact embedding $(T^2 \times I, \xi_{\pi}) \hookrightarrow (\Sigma
\times S^1, \xi)$, with $T^2 \times \{0\}$ isotopic to a boundary
component, where $\xi_{\pi}$ is a tight contact structure with convex
boundary and twisting $\pi$ (see~\cite[{\S}~2.2.1]{H1} for the
definition of twisting). By the results of~\cite[Section~2]{GLS2}, we
are free to change the restriction $\xi\vert_{M-(U_1\cup U_2\cup
U_3)}$ without changing the isotopy class of $\xi$, as long as we
replace $\xi\vert_{M-(U_1\cup U_2\cup U_3)}$ by a tight, appropriate
contact structure with minimal convex boundaries having the same
infinite boundary slopes.

We claim that this contact structure can be chosen to be the contact
complement $(C,\eta)$ of standard neighborhoods of the three
Legendrian unknots $L_1, L_2, L_3$ in Figure~\ref{f:structures}. In
fact, the contact structure obtained by performing contact
$(-1)$--surgeries along the knots $L_1,L_2, L_3$ is tight, due to the
cancellation of $(\pm 1)$--surgeries of opposite signs along
Legendrian push--offs~\cite{DG2}. Therefore $\eta$ is tight, and it is
easy to check that $\eta$ is also appropriate because it extends to a
tight contact structure on a closed 3--manifold obtained by filling
the neighborhoods of $L_1,L_2$ and $L_3$. It is obvious that the
boundary components of $(C,\eta)$ are minimal and convex. To check
that the boundary components have infinite boundary slopes, it is
enough to observe that there is a product structure $C\cong \Si\x S^1$
such that (i) a fiber $F$ of the projection $C\to\Si$ is Legendrian
and has twisting number zero, and (ii) all the $L_i$'s are Legendrian
pushoffs of $F$.
\end{proof}

\sh{Open book decompositions}

According to a recent result of Giroux~\cite{Gi}, isotopy classes of
contact structures are in one-to-one correspondence with suitable
equivalence classes of open book decompositions of the underlying
3--manifold, cf. also~\cite{Et}.  Recall that an open book
decomposition $(Y,f)$ of a closed 3--manifold $Y$ amounts to a link
$L\subset Y$ and a fibration $f\colon Y-L$ such that the closure of
every fiber $f^{-1}(t)$ (called a \emph{page}) provides a Seifert
surface for $L$. The open book decomposition is called \emph{planar}
if the genus of the page is zero. We will also call \emph{planar} any
contact structure compatible with a planar open book decomposition.
The significance of this notion is evident from the following result.

\begin{thm}[Abbas--Cieliebak--Hofer~\cite{ACH}]\label{t:hofer}
If $\xi $ is a contact structure compatible with a planar open book 
decomposition then it satisfies the Weinstein conjecture, that is,
any Reeb vector field of $\xi $ admits a periodic orbit. \qed
\end{thm}

Necessary conditions for a contact structure to be planar were found
in~\cite{Et1, OSS}. We will prove Theorem~\ref{t:planob} using
Proposition~\ref{p:model}. Before dwelving into the proof, we describe
some consequences of Theorem~\ref{t:planob}, some of which immediately
imply Corollaries~\ref{c:planar=L},~\ref{c:allplanar}
and~\ref{c:weinstein}. First of all, we have the following

\begin{cor}\label{c:planar-first}
If $M=M(e_0;r_1,r_2,r_3)$, $e_0\geq -1$ and $M$ is an $L$--space, then
each contact structure on $M$ is planar and therefore satisfies
Weinstein's conjecture.
\end{cor}

\begin{proof}
It is known that overtwisted contact structures are planar~\cite{Et1},
hence we can focus on tight structures only. For small Seifert fibered
rational homology 3--spheres with $e_0\geq 0$ it was shown by
Wu~\cite{Wu0} that each tight contact structure has maximal twisting
equal to zero, therefore Theorem~\ref{t:planob} shows that on those
manifolds every contact structure is planar.

Suppose now that $M$ is an $L$--space and $e_0(M)=-1$. By
Theorem~\ref{t:equiv} $M$ admits no transverse contact structures,
hence Theorem~\ref{t:ghig} shows that every tight contact structure on
$M$ has maximal twisting equal to zero.  Theorem~\ref{t:planob}
therefore implies that every tight contact structure on $M$ is planar.
\end{proof}

\begin{rem}\label{r:poincare}
In~\cite{Sc} Stephan Sch\"onenberger proved that if
$e_0\leq -3$ then every contact structure on $M=M(e_0;r_1,r_2,r_3)$ is
planar. Therefore, since in this case $M$ is always an $L$--space, the
statement of Corollary~\ref{c:planar-first} holds if $e_0\leq -3$.  On
the other hand, the Poincar\'e sphere admits a Seifert fibration with
$e_0=-2$, is an $L$--space, but its unique Stein fillable contact
structure is not planar because the intersection form of a filling
is not diagonalizable~\cite{Et1,OSS}.
\end{rem}

In view of Corollary~\ref{c:planar-first}, Theorem~\ref{t:equiv}
implies the following characterization of $L$--spaces of the form
$M(e_0;r_1,r_2,r_3)$.

\begin{cor}\label{c:planar-second}
The 3--manifold $M=M(e_0;r_1,r_2,r_3)$ is an $L$--space if and only
if one of the following holds:
\begin{itemize}
\item
each contact structure on $M$ is planar;
\item
each contact structure on $-M$ is planar.
\end{itemize}
\end{cor}

\begin{proof}
If $e_0(M)\geq 0$ then by Corollary~\ref{c:planar-first} every contact
structure on $M$ is planar.  On the other hand, it is known that if
$e_0(M)\geq 0$ or $e_0(M)\leq -3$ then $M$ is an $L$--space. Since
$e_0(-M)=-3-e_0(M)$, this immediately implies the statement if
$e_0(M)\geq 0$ or $e_0(M)\leq -3$. Suppose now that $M$ is an
$L$--space and, up to changing its orientation, $e_0(M)=-1$. Then, by
Corollary~\ref{c:planar-first} every contact structure on $M$ is
planar.

To finish the proof, we may assume without loss that $e_0(M)=-1$ and
$M$ is not an $L$--space. According to Theorem~\ref{t:equiv} the
manifold $M$ admits a taut foliation, which gives rise to a contact
structure $\xi$ having a symplectic semi--filling. By~\cite{El, Et2}
one can use a symplectic cap to construct a symplectic filling of
$(M,\xi)$ with $b_2^+>0$. Then, according to~\cite{Et1, OSS} $\xi$ is
not planar.
\end{proof}

Let $T_{p,np+1}$ denote the positive $(p,np+1)$--torus knot, and let
\[
M_{p,pn+1}=-S^3_{p^2n-pn-1}(T_{p,pn+1}). 
\]
Proposition~\ref{p:nonplanar} below shows that the non--planarity of
the Stein fillable contact structure on the Poincar\'e sphere
(cf. Remark~\ref{r:poincare}) is a non--isolated phenomenon.

\begin{prop}\label{p:nonplanar}
For $p\geq 2$ and $n\geq 1$ the rational homology 3--sphere
$M_{p,np+1}$ is an $L$--space and it
carries non--planar contact structures.
\end{prop}

\begin{proof}
Since the slice genus of $T_{p,np+1}$ is $p^2n-pn-1$, the fact that
$M_{p,pn+1}$ is an $L$--space follows from~\cite[Proposition~4.1]{LSgt}.

Let $W_{p,n}$ denote the canonical plumbing 4--manifold associated to
$M_{p,np+1}$ as in Section~\ref{s:sec}. Since all the weights of the
plumbing are $\leq -2$ (and all knots are unknots), $W_{p,n}$ supports
Stein structures inducing tight contact structures on $M_{p,n}$. The
proof of~\cite[Theorem~4.1]{Et1} shows that if any of these structures
is planar, then $W_{p,n}$ smoothly embeds in a connected sum of
$\overline\CP^2$'s. But the argument given in the proof
of~\cite[Proposition~4.1]{LS2} shows that the intersection form
$Q_{W_{p,n}}$ does not embed into a diagonal lattice. Therefore, none
of the contact structures filled by $W_{p,n}$ are planar.
\end{proof}

\begin{proof}[Proof of Theorem~\ref{t:planob}]
Let $\xi$ be a contact structure with twisting number equal to zero on
$M(-1; r_1,r_2,r_3)$. By Proposition~\ref{p:model} this structure is
obtained by performing a (possibly rational) contact surgery along the
five--component Legendrian link $\cal L$ of Figure~\ref{f:special}
(for $k=3$). According to the algorithm outlined in
Section~\ref{s:sec}, $\xi$ is obtained by contact $(\pm 1)$--surgery
on a Legendrian link $\tilde{\cal L}$ obtained from $\cal L$ by
successively taking pushoffs and Legendrian stabilizations of (some
of) its components. It is well known that performing contact $(\pm
1)$--surgery on a Legendrian knot which sits on a page of a compatible
open book with contact framing equal to the page framing yields an
open book of the same genus compatible with the resulting contact
structure. Therefore, it suffices to show that $\tilde{\cal L}$ sits
on a union of pages of a planar open book for $S^3$ compatible with
the standard contact structure. This can be proved by an argument very
similar to the one used in~\cite{Sc} to prove that each contact
structure on a lens space is planar. From now on, we refer
to~\cite{Et} for standard facts on contact structures and their
compatible open books. Start with the open book decomposition of $S^3$
compatible with the standard contact structure whose page is an
annulus and whose binding is a Hopf link.  By applying the Legendrian
realization principle, a Legendrian unknot with maximal
Thurston--Bennequin invariant together with four of its pushoffs can
be realized on five distinct pages of this open book, so that contact
framings and page framings coincide. Up to positively stabilizing the
open book in the sense of Giroux, we can realize in the same way any
Legendrian stabilizations of these unknots on distinct pages of a
planar open book compatible with the standard contact structure on
$S^3$. Any number of pushoffs of the stabilized knots can then be
realized on distinct pages, any further Legendrian stabilization can
be realized on further planar stabilizations of the open book, and so
on. This construction clearly establishes what we need.
\end{proof}


\begin{thebibliography}[AAA]

\bibitem{ACH}
{\bf  C. Abbas, K. Cieliebak and H. Hofer},
{\it The Weinstein conjecture for planar contact structures in
dimension three},
arxiv:math.SG/0409355.

\bibitem{Bu}
{\bf G. Burde and  H. Zieschang},
{\em Knots},
de Gruyter Studies in Mathematics, {\bf5}.
Walter de Gruyter \& Co., Berlin, 1985.

\bibitem{DG2}
{\bf F. Ding and H. Geiges},
{\em A Legendrian surgery presentation of contact 3-manifolds},
Math. Proc. Cambridge Philos. Soc. {\bf136} (2004) 583--598.

\bibitem{DGS}
{\bf F. Ding, H. Geiges and A. Stipsicz},
{\it Surgery diagrams for contact 3--manifolds}, Turkish J. Math. 
{\bf28} (2004) 41--74.

\bibitem{El}
{\bf Y. Eliashberg},
{\it A few remarks about symplectic filling},
Geom. Topol. {\bf 8} (2004) 277--293.

\bibitem{ET} 
{\bf Y. Eliashberg and W. Thurston},
{\it Confoliations},
University Lecture Series, {\bf 13}, 
American Mathematical Society, Providence, RI, 1998.

\bibitem{Et}
{\bf J. Etnyre},
{\it Lectures on open book decompositions and contact structures},
arxiv:math.SG/0409402.

\bibitem{Et1}
{\bf J. Etnyre},
{\it Planar open book decompositions and contact structures},
Int. Math. Res. Not. {\bf 2004}, no. 79, 4255-4267.

\bibitem{Et2}
{\bf J. Etnyre},
{\it On symplectic fillings},
Algebr. Geom. Topol. {\bf 4} (2004) 73--80.

\bibitem{EH1}
{\bf J. Etnyre and K. Honda},
{\it Tight contact structures with no symplectic fillings},  
Invent. Math. {\bf148} (2002), 609--626.

\bibitem{EH2}
{\bf J. Etnyre and K. Honda},
{\it On the nonexistence of tight contact structures},
Ann. of Math. {\bf153} (2001), 749--766.

\bibitem{Gh}
{\bf P. Ghiggini},
{\it Tight contact structures with negative maximal twisting on
small Seifert manifolds}, in preparation.

\bibitem{Gi} 
{\bf E. Giroux} 
{\it G\'eom\'etrie de contact: de la dimension trois vers les
dimensions sup\'erieures},
Proceedings of the International Congress of Mathematicians, Vol. II
(Beijing, 2002), 405--414, Higher Ed. Press, Beijing, 2002.

\bibitem{GLS}
{\bf P. Ghiggini, P. Lisca and A. Stipsicz},
{\it Classification of tight contact structures on
small Seifert 3--manifolds with $e_0\geq 0$},
arxiv:math.SG/0406080, to appear in Proc. Amer. Math. Soc.

\bibitem{GLS2}
{\bf P. Ghiggini, P. Lisca and A. Stipsicz},
{\it Tight contact structures on some small Seifert fibered 3--manifolds},
in preparation

\bibitem{Go}
{\bf R. Gompf},
{\it Handlebody constructions of Stein surfaces},
Ann. of Math. {\bf148} (1998), 619--693.

\bibitem{H1}
{\bf K. Honda},
{\it On the classification of tight contact structures, I.},
Geom. Topol. {\bf4} (2000) 309--368.

\bibitem{H2}
{\bf K.~Honda}, 
{\it Confoliations transverse to vector fields}, preprint, 
preliminary version, http://math.usc.edu /~khonda/research.html

\bibitem{JN} 
{\bf M.~Jankins and W.~Neumann}, 
{\it Rotation Numbers of Products of Circle Homeomorphisms}, 
Math. Ann. {\bf 271} (1985) 381--400.

\bibitem{LM}
{\bf P. Lisca and G. Mati\'c},
{\it Transverse contact structures on Seifert 3-manifolds},
 Algebr. Geom. Topol. {\bf4}  (2004), 1125--1144.

\bibitem{LS1}
{\bf P. Lisca and A. Stipsicz},
{\it An infinite family of tight, not semi--fillable contact 3-manifolds}, 
Geom. Topol. {\bf7} (2003), 1055--1073.

\bibitem{LS1.5}
{\bf P. Lisca and A. Stipsicz}, 
{\it Seifert fibered contact three--manifolds via surgery}, 
Alg. Geom. Topol. {\bf 4} (2004), 199--217.

\bibitem{LSgt}
{\bf P. Lisca and A. Stipsicz}, 
{\it Ozsv\'ath--Szab\'o invariants and tight contact three-manifolds. I},
Geom. Topol.  {\bf8}  (2004), 925--945.

\bibitem{LS2}
{\bf P. Lisca and A. Stipsicz},
{\it Ozsv\'ath--Szab\'o invariants and tight contact three-manifolds. II},
arxiv:math.SG/0404136.

\bibitem{McW} 
{\bf J.~McCarthy and J.~Wolfson}, 
{\it Symplectic gluing along hypersurfaces and resolution 
of isolated orbifold singularities}, 
Invent. Math. {\bf 119} (1995), 129--154.

\bibitem{Na} 
{\bf R.~Naimi}, 
{\it Foliations transverse to fibers of Seifert manifolds},
Comment. Math. Helv. {\bf 69} (1994), 155--162.

\bibitem{Ne}
{\bf A. N\'emethi},
{\it On the Ozsv\'ath--Szab\'o invariant of negative definite 
plumbed 3-manifolds}, 
Geom. Topol. {\bf9} (2005), 991--1042.

\bibitem{NR} 
{\bf W. Neumann and F. Raymond},
{\it Seifert manifolds, plumbing, $\mu $--invariant and orientation 
reversing maps}, 
Algebraic and geometric topology (Proc. Sympos., Univ. California, 
Santa Barbara, Calif., 1977), pp. 163--196, 
Lecture Notes in Math. {\bf 664}, Springer, Berlin, 1978.

\bibitem{Or}
{\bf P. Orlik}, 
{\it Seifert manifolds}, 
Lecture Notes in Math., Vol. {\bf 291}. 
Springer-Verlag, Berlin-New York, 1972. 

\bibitem{OSS}
{\bf  P. Ozsv\'ath, A. Stipsicz and  Z. Szab\'o},
{\it Planar open books and Floer homology},
arxiv:math.SG/0504403.

\bibitem{OSzF1} {\bf P. Ozsv\'ath and Z. Szab\'o}, {\it Holomorphic
disks and topological invariants for closed three-manifolds}, 
Ann. of Math. {\bf159}  (2004), 1027--1158.

\bibitem{OSzF2} {\bf P. Ozsv\'ath and Z. Szab\'o}, {\it Holomorphic
 disks and three--manifold invariants: properties and applications}, 
Ann. of Math. {\bf159} (2004), 1159--1245.

\bibitem{abs}
{\bf P. Ozsv\'ath and Z. Szab\'o},
{\it Absolutely graded Floer homologies and intersection forms
for four--manifolds with boundary},
Adv. Math. {\bf173} (2003), 179--261.

\bibitem{OSzF4}
{\bf P. Ozsv\'ath and Z. Szab\'o},
{\it Holomorphic triangles and invariants of smooth $4$--manifolds},
arXiv:math.SG/0110169. 

\bibitem{OSzplum}
{\bf P. Ozsv\'ath and Z. Szab\'o},
{\it On the Floer homology of plumbed three-manifolds},
Geom. Topol. {\bf7} (2003), 185--224.

\bibitem{OSzcont}
{\bf P. Ozsv\'ath and Z. Szab\'o},
{\it Heegaard Floer homologies and contact structures},
arxiv:math.SG/0210127.

\bibitem{OSzlens}
{\bf P. Ozsv\'ath and Z. Szab\'o},
{\it On knot Floer homology and lens space surgery},
arxiv:math.GT/0303017.

\bibitem{OSzsurvey} {\bf P. Ozsv\'ath and Z. Szab\'o}, {\it Heegaard
diagrams and holomorphic disks}, Different faces of geometry,
301--348, Int. Math. Ser. (N. Y.), Kluwer/Plenum, New York, 2004.

\bibitem{OSzgen}
{\bf P. Ozsv\'ath and Z. Szab\'o},
{\it Holomorphic disks and genus bounds},
Geom. Topol. {\bf8} (2004), 311--334.

\bibitem{OSrat}
{\bf P. Ozsv\'ath and Z. Szab\'o},
{\it Knot Floer homology and rational surgeries},
arxiv:math.GT/0504404.

\bibitem{Sc}
{\bf S. Sch\"onenberger},
{\it Planar open books and symplectic fillings},
Ph.D. thesis, University of Pennsylvania (2005).

\bibitem{Wu0}
{\bf H. Wu},
{\it Legendrian Vertical Circles in Small Seifert Spaces},
arxiv:math.GT/0310034.

\bibitem{Wu}
{\bf H. Wu},
{\it Tight Contact Small Seifert Spaces with $e_0\neq0,-1,-2$},\newline
arxiv:math.GT/0402167.

\end{thebibliography}
\end{document}